\documentclass[11pt, twoside]{article}
\usepackage{amsmath,amsthm}
\usepackage{amsfonts}
\usepackage{amssymb,latexsym}
\usepackage{enumerate}
\newtheorem{theorem}{Theorem}[section]
\newtheorem{remark}[theorem]{Remark}
\newtheorem{proposition}[theorem]{Proposition}
\newtheorem{lemma}[theorem]{Lemma}
\newtheorem{definition}[theorem]{Definition}

\newtheorem{corollary}[theorem]{Corollary}
\numberwithin{equation}{section}
\def\ind{1{\hskip -3 pt}\hbox{\textsc{I}}}
\def\n{\noindent}

\def\Ga{\Gamma}
\def\de{\delta}

\def\o{\omega}
\def\O{\Omega}

\def\bR{\mathbb R}
\def\al{\alpha}

\def\w{\wedge}
\def\cn{\mathbb C^n}
\def\om{\omega}
\def\pa{\partial}

\def\R{\mathbb{R}}
\def\bC{\mathbb C}

\def\ga{\gamma}
\def\wed{\wedge}
\def\Om{\Omega}

\def\d{\partial}

\begin{document}
	\setlength{\baselineskip}{18truept}
	\pagestyle{myheadings}
	
	\title {Bounded solution to Hessian type equations for $(\o,m)-\beta$-subharmonic functions on a ball in $\mathbb C^n$}
	\author{
	 Hoang Thieu Anh*,	Le Mau Hai**, Nguyen Quang Dieu** and Nguyen Van Phu*** \\
		*Faculty of Basic Sciences, University of Transport and Communications, 3 Cau Giay,\\ Dong Da, Hanoi, Vietnam.\\
		 **Department of Mathematics, Hanoi National University of Education,\\ Hanoi, Vietnam.\\
		*** Faculty of Natural Sciences, Electric Power University,\\ Hanoi,Vietnam;
		\\E-mail: anhht@utc.edu.vn, mauhai@hnue.edu.vn \\ngquang.dieu@hnue.edu.vn and  phunv@epu.edu.vn 
	}

	\date{}
	\maketitle
	
	\renewcommand{\thefootnote}{}
	
	
	
	\renewcommand{\thefootnote}{\arabic{footnote}}
	\setcounter{footnote}{0}

\begin{abstract}
	\n
In this paper, we study Hessian type equations for $(\o,m)-\beta$-subharmonic functions on a ball in $\mathbb{C}^n$, where $\beta=dd^c\|z\|^2=\frac{i}{2}\sum\limits_{j=1}^n dz_j\w d\bar{z}_j$ is the flat metric on $\cn$.
Using the  recent results in  \cite{KN23b}, we are able to show
the existence of bounded solutions for Hessian type equations. 
	\end{abstract}
	
	\section{Introduction}
	
Let $\Om\subset \bC^n$ be a ball and $\om$ be a smooth real $(1,1)$-form on $\mathbb C^n$.  
Let $1\leq m\leq n$ be an integer and consider  a function $u\in C^2(\Om,\bR)$. The complex Hessian operator  with respect to $\om$ acts on  $u$ by
$$
H_{m,\o}(u) = (\o+dd^c u)^m \wedge \beta^{n-m},
$$
where $\beta=dd^c||z||^2$ is a standard K\"ahler form in $\mathbb C^n$.

For $\om = 0,$    B\l ocki \cite{Bl05}  defined the Hessian operator acting on bounded $m$-subharmonic functions (not necessarily smooth).
In this case,
the Hessian operator $H_{m,0}(\bullet)$ are positive Radon measures which is stable under monotone sequences, and that the homogeneous Dirichlet problem is solvable on a ball in $\cn$. To make it convenient for presentation, from now on we use symbol $H_m(\bullet)$ instead of $H_{m,0}(\bullet)$.

According to \cite{GN18} and \cite{KN23c}, a function $u:\Om \to [-\infty,+\infty)$ is called $(\o,m)-\beta$-subharmonic ($(\o,m)-\beta$-sh for short) if $u\in L^1_{\rm loc}(\Om,\beta^n)$ and for any collection $\ga_1,...,\ga_{m-1} \in \Ga_m(\beta)$, 
$$
(\o+dd^c u )\wed \ga_1 \wed \cdots \wed \ga_{m-1} \wed \beta^{n-m}\geq 0
\eqno(1.1)$$
in the sense of currents, where the positive cone $\Ga_m(\beta)$, associated to $\beta$, is defined as 
$$
\Ga_m(\beta):= \{\ga: \, \ga\ \ \text{is a (1,1)}-\text{real form}, \ga^k \wed \beta^{n-k} (z)>0  \quad\text{for } k=1,...,m\}.
$$
The cone of $(\o,m)-\beta$-sh (resp. negative) functions on $\Om$ is denoted by $SH_{ \om,m} (\Om$) (resp. $SH_{ \om,m}^{-} (\Om)$). In the case when $\o=0,$ to make it convenient for readers, from now on we use symbol $SH_{m} (\Om$) (resp. $SH_{m}^{-} (\Om)$) instead of $SH_{ 0,m} (\Om$) (resp. $SH_{ 0,m}^{-} (\Om)$).\\
For a $C^2$ function $u$, the inequality (1.1) is equivalent to inequalities

\begin{equation}
	(\o+dd^cu)^k\w\beta^{n-k}\geq 0\,\,\text{for}\,\,k=1,\ldots,m	
	\end{equation}
Let $\rho$ be a strictly plurisubharmonic function on $\O$ such that $dd^c\rho\geq \o$ on $\Omega,$ then $u+\rho$ is a $m$-sh function on $\O$. We write $\tau=dd^c\rho-\o$ which is a smooth $(1,1)$-form, then $\o+dd^cu=dd^c(u+\rho)-\tau.$ We define 
\begin{equation}\label{e1.3}(\o+dd^cu)^m\w\beta^{n-m}:=\sum_{k=0}^m\binom{m}{k}(-1)^{m-k}[dd^c(u+\rho)^k]\w\tau^{m-k}\w\beta^{n-m}.
	\end{equation}
Following the tradition inductive method as in \cite{BT1}, the Hessian operator $H_{\o,m} (u)$ can be defined
over the class of locally bounded $(\o,m)-\beta$-sh functions $u$ as  positive Radon measure which put no mass on $m$-polar sets. See Section 9 in \cite{KN23c} and \cite{GN18} for details. Moreover, in \cite{GN18}, Dongwei Gu and  Ngoc Cuong  Nguyen developed a pluripotential theory analogous to the classical work \cite{BT1} where fundamental properties like continuity under monotone sequences of the Hessian operator, quasi-continuity and the comparison principle,... are also true. Note that Remark 3.19 (b) in \cite{GN18} implies that the class  every bounded $(\o,m)-\beta$-sh functions can be approximated by a decreasing sequence of smooth $(\o,m)-\beta$-sh functions.

\n 
From now on, unless otherwise specified, we always assume that 
$\mu$ is a positive Radon measure on 
$\Om ,\phi:\partial \O \to \mathbb{R}$ is continuous function 
and $F(t,z): \R\times \Om \to [0,+\infty)$ is a continuous function
such that there exists a function $G \in L_{loc}^1 (\Om, \mu)$ satisfying
\begin{equation}\label{condi}F(t,z) \le G(z), \ \forall (t, z) \in \R \times \Om.
	\end{equation}
	
The main goal of this paper is to seek for a solution $u$ of the following Hessian type equation
	\begin{equation}\tag{*}
		\begin{cases}
			u\in SH_{\om,m}(\Om)\cap L^{\infty}(\overline{\Om})\\
			H_{m,\o}(u)=F(u,z)d\mu\\
			\lim\limits_{B\ni z\to x}u(z)=\phi (x), \forall x\in\pa \Om. 
		\end{cases}
	\end{equation}
In the case when $m=n,$ this equation is called Monge-Amp\`ere type equation. This equation attracts the interest of many mathematicians.
In the case when $ \omega=0,$  $d\mu=dV_{2n}$, Bedford and Taylor \cite{BT79} proved that there exist a unique solution to the problem (*) with $F(t,z)\in C^{0}(\R\times\overline{\Om})$ such that $F^{\frac{1}{n}}$ is convex and non-decreasing in $t$. Afterthat, Cegrell \cite{Ce84} extended this result by showing that (*) has solution when $F(t,z)$ is bounded function which is continuous with respect to the first variable for every fixed $z\in\O.$ When $ \o=0$ and the measure $\mu$ is more general, S. Ko{\l}odziej \cite{K00} showed that if $d\mu=(dd^cv)^n$, where $v\in PSH(\Omega)\cap L^{\infty}(\Omega)$ and $\lim\limits_{z\to x}v(z)=\phi(x)$ for $x\in\partial \O$ as well as $F(t,z) $ is a bounded function which is non-decreasing and continuous in the first variable, $d\mu$-measurable in the second one then the problem (*) has a unique solution. Recently, in the case when $\o$ is more general, the problem (*) has been discussed in \cite[Theorem 2.3]{KN23b}. Note that in the above paper, Kolodziej and Cuong consider the function $F(t,z)$ to be bounded, non-decreasing with respect to the first variable, and they use a familiar tool to prove the existence of a solution: the Schauder fixed-point Theorem.\\
Motivated by the above results, in this paper, we want to extend the results of Kolodziej and Cuong from Monge-Amp\`ere-type equations to Hessian type equations. Our novelty is that we do {\it not} assume monotonicity of $t \mapsto F(t,z)$ as in \cite{KN23b} and the hypothesis $F(t,z)$ is a bounded function in \cite{KN23b} is replaced by the hypothesis $F(t,z)$ is bounded by a local integrable function with respect to the measure $\mu$ as in inquality \eqref{condi}  which is more general. Moreover, by using convergence in $(\o,m)$-capacity, we can avoid the use of Schauder  fixed point Theorem. To the best of our knowledge, this approach seems to be new.

Our first main result is the proof of the following Theorem:
\begin{theorem} \label{th1.1}
	Assume that the following conditions hold true:

\n 
(a) There exists $v\in SH_{m}^-(\Om)\cap L^{\infty}(\Om)$ satisfying
$$\lim\limits_{\Om \ni z\to x}v(z)=0\,\forall x\in\pa \Om \ \text{and}\ G\mu\leq H_{m}(v);$$ 
\n	
(b) For every $m-$polar subset $E$ of $\Om$ we have $\mu (E \cap \{G=0\})=0.$

Then the problem (*) has a solution.\\
 Moreover, if $ F(t,z)$ is a non-decreasing function with respect to the first variable for every 
	$z \in \O \setminus X$ where $X$ is a Borel set with $C_m (X)=0,$ then the problem (*) has an unique solution.
\end{theorem}

\n
Finally we also concern  with a sort of stability of solutions of (*). This result should be compared with Theorem 2.7 in \cite{KN23b}, where the authors establised a result on stability of solutions to (*) in the context of 
$\om-$plurisubharmonic functions on a compact Hermitian manifolds with boundary $(M,\omega)$.
\begin{theorem} \label{th1.2}
Let $F, F_j: \R \times \Om \to [0, \infty) \ (j \ge 1)$ be a set
of continuous functions that satisfies inequality \eqref{condi} and
that for every $z \in \Om \setminus Y,$ the sequence
$F_j (t, z)$ converges locally uniformly to $F (t,z),$ where $Y\subset\O$ is a $m$-polar set.
Moreover, suppose that the conditions (a) and (b) in Theorem \ref{th1.1} hold true.
For each $j,$ let $u_j$ be a solution of the equation
	\begin{equation}
\begin{cases}
u_j\in SH_{\om,m}(\Om)\cap L^{\infty}(\Om)\\
H_{m,\o}(u_j)=F_j(u_j,z)d\mu\\
\lim\limits_{\O\ni z\to x}u_j(z)=\phi (x), \forall x\in\pa \Om. 
\end{cases}
\end{equation}
Then there exists a subsequence of $\{u_j\}$ that converges in $(\o,m)$-capacity to a solution $u$ of the problem (*).
\end{theorem}
\n The paper is organized as follows. Besides the introduction, the paper has other two sections. In Section 2, following seminal works \cite{KN23c}, \cite{GN18},... we collect basic features of $(\o,m)-\beta$-sh functions on $\O.$
Most notably is the comparision principle for the Hessian operator $H_{m,\o} (u),$ we also recall the subsolution Theorem, a powerful tool to check existence of solution of Dirichlet problem for Hessian operator. Another important ingredient is the relationship between the convergence in $(\o,m)$-capacity and the convergence of corresponding Hessian operators.	This fact will allow us to bypass the use of Schauder
fixed point Theorem in the proof of Theorem \ref{th1.1}.
In Section 3, we supply in details the proofs of our main results.

\n 
{\bf Acknowledgments} 
The third name is supported by Grant number 101.02-2023.11 from the NAFOSTED (National Foundation for Science and Technology developements) program.\\
This work is written in our visit in VIASM in the Spring of 2025. We also thank VIASM for financial support and hospitality.

	\section{Preliminaries} 
Following Proposition 2.6 in  \cite{GN18}, we include below some basic properties of $(\o,m)-\beta$-sh functions. 
\begin{proposition} \label{prop:closure-max}
	\begin{enumerate}
		\item[(a)]
		$SH_{\o,n}(\O)\subset SH_{\o,n-1}(\O)\subset\cdots\subset SH_{\o,1}(\O)$
	\item[(b)]	If $u_1 \geq u_2 \geq  \cdots$ is a decreasing sequence of $(\o,m)-\beta$-sh functions, then $u := \lim_{j\to \infty} u_j$ is either $(\o,m)-\beta$-sh or $\equiv -\infty$.
		\item[(c)] 
		If $u, v$ belong to $SH_{ \om,m}(\Omega)$, then so does $\max\{u,v\}$.
		\item[(d)] 
		If $u, v$ belong to $SH_{m, \om}(\Omega)$ and satisfies $u \le v$ a.e. with respect to Lebesgue measure then $u \le v$ on $\Om.$
	\end{enumerate}
\end{proposition}
\n
Notice that $(d)$ follows from Lemma 9.6 and Definition 2.4 in \cite{GN18}  where $u+\rho$ and $v+\rho$ are viewed as $\al-$subharmonic function with any $(1,1)$ form $\al$ such that 
$
\al^{n-1} =\ga_1\w\cdots\w \ga_{m-1} \wed \beta^{n-m},$ where $\ga_1,\ldots,\ga_{m-1}$ is a certain form belonging to $\Ga_m (\beta)$.

\n We next recall the notion of capacity associated with Hessian operators of bounded $(\o,m)-\beta$-subharmonic functions (see equality (3.6) in \cite{GN18}).
	\begin{definition}
		For a Borel set $E\subset \Om,$ we set
		$$Cap(E):=\sup\left\{\int_{E}(\o+dd^cv)^m \wedge \beta^{n-m}: v\in SH_{\om,m}(\Om), 0\leq v\leq 1\right\},$$
		$$C_m(E):=\sup\left\{\int_{E}(dd^cv)^m \wedge \beta^{n-m}: v\in SH_{m}(\Om), 0\leq v\leq 1\right\}.$$
	\end{definition}
\n According to Lemma 3.5 in \cite{GN18}, there exists a constant $C$ depending on $\o$ such that \begin{equation}\label{e2.1}\frac{1}{C}Cap(E)\leq C_m(E)\leq C. Cap(E)\end{equation}
\begin{remark}\label{remark}{\rm
	Let $\{u_\alpha\}_{\alpha \in I} \subset SH_{\om,m}(\Omega)$ be a family which is locally uniformly bounded from above. Put $u(z) := \sup_\alpha u_\alpha(z)$. According to Proposition 2.6 in \cite{GN18}, the upper semicontinuous regularization $u^*$ is $(\om,m)-\beta$-sh. 
Moreover, using the same argument as Corollary 5.5 in \cite{KN23c}, we obtain that $Cap(E)=0$, where $E:=\{u<u^*\}.$ It follows from inequality \eqref{e2.1} that $C_m(E)=0.$ This means that $E$ is $m-$polar, i.e., there exists $v \in SH_{m} (\Om)$ such that $v \equiv -\infty$ on $E.$}
\end{remark}
\n Similar to Definition {\bf 4.2} in \cite{KN23c},  we have the following definition.
\begin{definition}
	A sequence of Borel functions $u_j$ in $\Omega$ is said to converge in $(\o,m)$-capacity (or in $Cap(.))$ to $u$ if for any $\delta>0$ and $K\Subset \Omega$ we have 
	$$\lim\limits_{j\to\infty}Cap(K\cap|u_j-u|\geq\delta)=0.$$
	\end{definition}
\begin{remark}\label{monocap} According to Corollary 4.11 in \cite{KN23c} and inequality \eqref{e2.1}, monotone convergence of locally uniformly bounded sequences of
$(\o,m)-\beta$-sh functions is convergence in $(\o,m)$-capacity.
\end{remark}

\n A major tool in pluripotential theory is the comparison principle.
We recalled a version of this result for bounded $(\o,m)-\beta$-sh. functions (see Corollary 3.11 in \cite{GN18}).

\begin{theorem}\label{comparison}
	Let $u,v$ be bounded $(\o,m)-\beta$-sh functions in  $\overline{\Omega}$ such that $\liminf\limits_{z\to\partial\Omega}(u-v)(z)\geq 0.$ Assume that $H_{m,\o}(v)\geq H_{m,\o}(u)$ in $\Omega.$ Then  $u\geq v$ on $\Omega.$
	\end{theorem}
\n 
Now, we will prove the following result. Note that this result was proved in the case of $m-$ subharmonic functions (see Lemma 4.5 in \cite{PD24}).
\begin{lemma}\label{bd1}
	Assume that $\mu$ vanishes on $m-$polar sets of $\Om$ and $\mu(\Omega) < \infty.$ Let $\{u_{j}\}\in SH_{\om,m}^{-} (\Omega)$ be a sequence satisfying the following conditions:
	
	\n 
	(i) $\sup\limits_{j \ge 1} \int\limits_{\Om} -u_jd\mu <\infty;$
	
	\n 
	(ii) $u_j \to u \in SH_{\om,m}^{-} (\Om)$ a.e. $dV_{2n}.$
	
	Then we have 
	$$\lim_{j \to \infty} \int\limits_{\Om} \vert u_j- u \vert d\mu=0.$$
	In particular $u_j \to u$ a.e. $d\mu$ on $\Om.$
\end{lemma}		
\begin{proof}
	Let $\rho\in C^2(\overline{\O})\cap PSH(\O)$  such that $dd^c\rho\geq \o$. Then $u_j+\rho$ and $u+\rho$ are $m$-subharmonic. Applying Lemma 4.5 in \cite{PD24} with functions $u_j+\rho$ and $u+\rho$ we get the desired.
	\end{proof}

\n According to Theorem 3.4 in \cite{GN18}, we know that $H_{m,\o} (u)$ is continuous
with respect to a monotone convergence of locally uniformly bounded sequences in $SH_{\om,m} (\Om).$
In the Proposition below, we present a more general result, where monotone convergence is replaced by convergence in $(\o,m)$-capacity.

\begin{proposition} \label{cap}
Let $u_j$ be a locally uniformly bounded sequence in $SH_{\om,m} (\Om).$ Assume that $u_j$ converges in $(\o,m)$-capacity
to  a locally bounded $u \in SH_{\om,m} (\Om)$. Then $H_{m,\o} (u_j)$ converges weakly to $H_{m,\o} (u).$
\end{proposition}
\begin{proof}
	
	Let $\rho$ be as in equality \eqref{e1.3}. We have $u_j+\rho$ and $u+\rho$ are $m$-subharmonic functions. It follows from the assumption $u_j$ converges in $(\o,m)$-capacity and inequality \eqref{e2.1} that $u_j+\rho$ is convergent in $m$-capacity to $u+\rho.$ By equality \eqref{e1.3} we only need to prove that $[dd^c(u_j+\rho)]^k\w\tau^{m-k}\w\beta^{n-m}$ converges weakly to $[dd^c(u+\rho)]^k\w\tau^{m-k}\w\beta^{n-m}$ for all $k=0,\cdots,m.$  Indeed, let $v$ be a strictly smooth plurisubharmonic function such that $dd^cv\geq \tau$ and $dd^cv\geq -\tau$ and let $\chi\in C^{\infty}_0(\Omega), 0\leq \chi\leq 1,$ we have 
\begin{align*}
	&\lim\limits_{j\to\infty}\bigg|\int_{\O}\chi\big[\big(dd^c(u_j+\rho)\big)^{k}-\big(dd^c(u+\rho)\big)^{k}\big]\w\tau^{m-k}\w\beta^{n-m}\bigg|\\
	&\leq\lim\limits_{j\to\infty}\bigg|\int\limits_{\{|u_j+\rho-(u+\rho)|>\varepsilon\}}\chi\big[\big(dd^c(u_j+\rho)\big)^{k}-\big(dd^c(u+\rho)\big)^{k}\big]\w\tau^{m-k}\w\beta^{n-m}\bigg|\\
	&+\lim\limits_{j\to\infty}\bigg|\int\limits_{\{|u_j+\rho-(u+\rho)|\leq\varepsilon\}}\chi\big[\big(dd^c(u_j+\rho)\big)^{k}-\big(dd^c(u+\rho)\big)^{k}\big]\w\tau^{m-k}\w\beta^{n-m}\bigg|.\\
	&\leq [\max_{\Omega}\chi]\lim\limits_{j\to\infty}\int\limits_{\{|u_j+\rho-(u+\rho)|>\varepsilon\}}\big[dd^c(u_j+\rho)\big]^{k}\w (dd^cv)^{m-k}\w\beta^{n-m}\\
	&+[\max_{\Omega}\chi]\lim\limits_{j\to\infty}\int\limits_{\{|u_j+\rho-(u+\rho)|>\varepsilon\}}\big[dd^c(u+\rho)\big]^{k}\w (dd^cv)^{m-k}\w\beta^{n-m}\\
	&+\bigg|\lim\limits_{j\to\infty}\int\limits_{\{|u_j+\rho-(u+\rho)|\leq\varepsilon\}}\chi\big[\big(dd^c(u_j+\rho)\big)^{k}-\big(dd^c(u+\rho)\big)^{k}\big]\w \tau^{m-k}\w\beta^{n-m}\bigg|.\\
	&\leq 2[\max_{\Omega}\chi] \lim\limits_{j\to\infty}C_m(\{|u_j+\rho-(u+\rho)|>\varepsilon\})\\
	&+\bigg|\lim\limits_{j\to\infty}\int\limits_{\{|u_j+\rho-(u+\rho)|\leq\varepsilon\}}\chi\big[\big(dd^c(u_j+\rho)\big)^{k}-\big(dd^c(u+\rho)\big)^{k}\big]\w (dd^cv)^{m-k}\w\beta^{n-m}\bigg|.\\
	&\leq \bigg|\lim\limits_{j\to\infty}\int\limits_{\{|u_j+\rho-(u+\rho)|\leq\varepsilon\}}\chi\big[\big(dd^c(u_j+\rho)\big)-\big(dd^c(u+\rho)\big)\big]\w T\w (dd^cv)^{m-k}\w\beta^{n-m}\bigg|.\\
	&\leq\lim\limits_{j\to\infty}\int_{\{|u_j+\rho-(u+\rho)|\leq\varepsilon\}}|u_j+\rho-(u+\rho)|dd^c\chi\w T\w (dd^cv)^{m-k}\w\beta^{n-m}.\\
	&\leq \varepsilon C_m(supp{\chi}).
\end{align*}
where $T=\sum_{l=0}^{k}\big(dd^c(u_j+\rho)\big)^l\w\big(dd^c(u+\rho)\big)^{k-l}.$\\
  
 The proof is complete.
\end{proof}

\n
Now we formulate the following subsolution Theorem  which plays a impotant role in our work (see Lemma 9.3 in \cite{KN23c}).
\begin{theorem}\label{sub}
Assume that there exists $v \in SH_{m} (\Om) \cap L^\infty (\Om)$ satisfying 
\[\label{eq:bounded-subsol}
H_{m}(v) \geq \mu, \quad \lim_{x\to \d \Om}  v(x) = 0.
\]	
Then, there exists a unique bounded $(\o,m)-\beta$-sh. function $u$ solving  
$$\lim_{z\to x} u(z) = \phi (x), \forall  x\in\d\Om,
	H_{m,\o}(u) = \mu\ \text{on } \Om.
	$$
\end{theorem}

	\section{ Weak solution to  Hessian type equation}

Firstly, we will prove the maximum principle in the class $SH_{\om,m}(\Om)\cap L^{\infty}(\Om).$

\begin{theorem}\label{thm 5.3}
	If $u, v \in SH_{\om,m}(\Om)\cap L^{\infty}(\Om)$, then 
	$$ \ind_{\{u>v\}} ( \omega + dd^c u)^m\w\beta^{n-m} = \ind_{\{u>v\}} ( \omega + dd^c \max(u,v))^m\w\beta^{n-m}. $$
\end{theorem}
\begin{proof}
	Fix $\rho \in \mathcal{P}_{m,\omega}(\Omega)$.  Write 
	$$H_{m,\o}(u)  = \sum_{k=0}^n \binom{m}{k} (-1)^{m-k} ( dd^c (u+\rho))^k \wedge (dd^c \rho - \omega)^{m-k}\w\beta^{n-m},$$
	$$H_{m,\o}( \max(u,v))  = \sum_{k=0}^n \binom{m}{k} (-1)^{m-k} ( dd^c ( \max(u,v)+\rho))^k \wedge (dd^c \rho - \omega)^{m-k}\w\beta^{n-m},$$
	
\n	Note that $(dd^c \rho - \omega)^{m-k}$ is a smooth $(m-k,m-k)$-form. Hence, according to Proposition 4.1 in \cite{S23} we can write
	$$  (dd^c \rho - \omega)^{m-k} = \sum\limits_{j\in J} f_j T_j,  $$
	where $J$ is a finete set, $(f_j)_{j\in J}$ are smooth functions with complex values  and $T_j=dd^cu^j_1\w\cdots\w dd^cu^j_{m-k}$  where, for every $j\in J$ and every $i=1,\cdots,m-k, u_i^j$ is a smooth negative plurisubharmonic function defined in a neighborhood of $\overline{\Omega}$. By linearity, it suffices to prove that
	$$\ind_{\{ u > v\}} ( dd^c (u+\rho))^k \wedge T_j\w\beta^{n-m}= \ind_{\{ u > v\}} (  dd^c \max(u+\rho,v+\rho))^k \wedge T_j \w\beta^{n-m}.  $$
	That is correct according to Lemma 4 in \cite{DE16}. The proof is complete.
\end{proof}

The following Corollary is an extension of Demailly's result. (see \cite[Proposition 6.11]{Dem89} in the case when $m=n$ and $\omega=0$).
\begin{corollary}\label{cor 5.5}
	Let $u, v \in SH_{\om,m}(\Om)\cap L^{\infty}(\Om)$. Then we have 
	$$ ( \omega + dd^c \max(u,v))^m\w\beta^{n-m} \geq \ind_{\{u> v \}} ( \omega + dd^c u)^m\w\beta^{n-m} + \ind_{\{ u \leq v \}} ( \omega + dd^c v)^m\w\beta^{n-m}. $$
\end{corollary}
\begin{proof}
	By Theorem \ref{thm 5.3}	we have 
	\begin{align*}
		&	( \omega + dd^c \max(u,v))^m\w\beta^{n-m} \\
		& \geq  \ind_{\{u>v\}} ( \omega + dd^c \max(u,v))^m\w\beta^{n-m} +  \ind_{\{u < v\}} ( \omega + dd^c \max(u,v))^m\w\beta^{n-m}\\
		& \geq  \ind_{\{u>v\}} ( \omega + dd^c u)^m\w\beta^{n-m} +  \ind_{\{u < v\}} ( \omega + dd^c v)^m\w\beta^{n-m}. 
	\end{align*}
	If $[( \omega + dd^c v)^m\w\beta^{n-m}] (\{u = v\}) = 0$, then the result follows. \\
In the case when $[( \omega + dd^c v)^m\w\beta^{n-m}] (\{u = v\}) \ne 0$,	since $(\omega + dd^c v)^m\w\beta^{n-m}$ 
	vanishes on $m$-polar sets,	the proof of \cite[Proposition 5.2]{HP17} shows that 
	$$ [(\omega + dd^c v)^m\w\beta^{n-m}](\{ u = v + t \}) = 0, \; \;  \forall t \in \mathbb{R}\setminus I_{\mu}, $$
	where $I_{\mu}$ is at most countable. Take $\varepsilon_j \in \mathbb{R}\setminus I_{\mu}$,  $\varepsilon_j \searrow 0$.   We have 
	
	$$[(\omega + dd^c v)^m\w\beta^{n-m}](\{ u = v + \varepsilon_j \}) = 0.$$
	This implies that
	
	\begin{align*} &( \omega + dd^c \max(u,v+ \varepsilon_j))^m\w\beta^{n-m} \\
		&\geq \ind_{\{u> v+ \varepsilon_j \}} ( \omega + dd^c u)^m\w\beta^{n-m} + \ind_{\{ u \leq v+ \varepsilon_j \}} ( \omega + dd^c v)^m\w\beta^{n-m}. 
	\end{align*}
	
	Let $\varepsilon_j \rightarrow 0$, according to Theorem 3.4 in \cite{GN18} and the Lebesgue Monotone Convergence Theorem, we get the desired.
\end{proof}

We need a version of the comparision principle. 
See  Proposition 2.2 in \cite{KN23b} for an analogous result for quasi-plurisubharmonic  functions.

	\begin{proposition}\label{md4}
Let $\nu \ge \mu$ be positive Radon measures on $\Om.$ 
Assume that	$ t \mapsto F(t,z)$ is a non-decreasing function in $t$ 
for all $z \in \Om \setminus Z,$ where $Z \subset \Om$ is a Borel set with $C_m (Z)=0.$
Let
		$u,v\in SH_{\om,m}(\Om)\cap L^{\infty}(\Om)$ be functions satisfying the following conditions: 

\n 
(i)
$\liminf\limits_{z\to\pa \Om}(u-v)(z)\geq 0;$

\n
(ii) $H_{m,\o}(u)=F(u,z)\mu, H_{m,\o}(v)=\tilde F(v,z)\nu,$ where $\tilde F \ge F$ is a
measurable function on $\Om.$
  
 Then $u\geq v$ on $\Om.$
	\end{proposition}

	\begin{proof}

		It follows from Corollary \ref{cor 5.5} that 
	\begin{align*}
		&(\omega + dd^c \max(u,v))^m\w\beta^{n-m} \\
		&\geq \ind_{\{u> v \}} ( \omega + dd^c u)^m\w\beta^{n-m} + \ind_{\{ u \leq v \}} ( \omega + dd^c v)^m\w\beta^{n-m} \\
		&= \ind_{\{u> v \}} F(u,.) d\mu  + \ind_{\{ u \leq v \}} \tilde{F}(v,.) d\nu \\
		&\geq \ind_{\{u> v \}} F(u,.) d\mu  + \ind_{\{ u \leq v \}} F(v,.) d\mu \\
		&= F(\max(u,v),.) d\mu \\
		&\geq F(u,.) d\mu = (\omega + dd^c u)^m\w\beta^{n-m},
	\end{align*}
	where the last inequality follows from the fact that the function $F$ is non-deceasing in the first variable. According to 
	Theorem \ref{comparison}, we infer that  $u \geq v$. The proof is complete.

	\end{proof}
\n	Now we give the proof of  Theorem \ref{th1.1}.

	\begin{proof} 
First, we show that $\mu$ puts no mass on $m-$polar subsets of $\Om.$ Indeed, for every  $m-$polar subset $E$ of $\Om,$
 by the assumption (a) we deduce that $(G\mu) (E)=0.$ By hypothesis (b), we obtain $\mu(E)=0.$
		
Next, according to Theorem \ref{sub}, there exists $h\in SH_{\om,m}(\Om)\cap L^{\infty}(\Om)$ such that
\begin{equation} \label{eqh}
H_{m,\o}(h)=0,\, h=\phi \ \text{on}\ \pa \Om.
\end{equation}
Set $u_0:= h$,	
 we have 
$$F(u_0, z)d\mu \le Gd\mu \le H_{m,\o} (v).$$
Thus, by Theorem \ref{sub}, we obtain aa unique $u_1 \in SH_{\om,m}(\Om)\cap L^{\infty}(\Om)$ satisfying
		 $$\lim\limits_{\Om \ni z\to x} u_1 (z)=\phi (x)\,\forall x\in\pa \Om, \
		H_{m,\o}(u_1)=F(u_0,z)d\mu.$$ 
		Since 
		$$H_{m,\o} (v+u_0) \ge  H_{m} (v) \ge Gd\mu \ge F(u_0,z)d\mu=H_{m,\o} (u_1),$$
	so using Theorem \ref{comparison} we obtain $v+u_0 \le u_1 \le u_0$ on $\Om.$
Continuing this process, we get a sequence $\{u_j\}	\in SH_{\om,m} (\Om) \cap L^\infty (\Om)$ that satisfies
	$$\lim\limits_{\Om \ni z\to x} u_{j+1} (z)=\phi (x)\,\forall x\in\pa \Om, \
	H_{m,\o}(u_{j+1})=F(u_j, z)d\mu.$$ 
We will show that $u_j$ converges in $(\o,m)$-capacity to a solution $\tilde{u}$ of $(*).$	
		
	 Since the sequence $\{u_j\}$ is uniformly bounded in $L^{\infty} (\Om),$ in view of Lemma 9.12 in \cite{GN18},
after swithching to a subsequence, we may achieve that
  $u_j \to u \in SH_{\om,m} (\Om)$ a.e. ($dV_{2n}$).		
On the other hand, by Lemma \ref{bd1}, we get $u_j \to u$ a.e. ($d\mu$).		
Indeed, it suffices to apply this lemma to relatively compact open subsets of $\Om$ and then applying
a diagonal process.
Now for $z \in \Om,$ we define the following sequences of non-negative measurable functions
	$$\varphi^1_j (z):= \inf_{k \ge j} F(u_k (z),z),
	\varphi^2_j (z):= \sup_{k \ge j} F(u_k (z),z).$$
	Then, by the inequality \eqref{condi}, we deduce that:\\
	\n 
	(i) $0 \le \varphi^1_j (z) \le  F(u_j (z),z) \le \varphi^2_j (z) \le G(z)$ for $j \ge 1;$
	
	\n 
	(ii) $\lim\limits_{j \to \infty} \varphi^1_j (z)=\lim\limits_{j \to \infty} \varphi^2_j (z)= F(u(z),z)$ a.e. ($d\mu$).
	
\n 	
Now, we may apply again Theorem \ref{sub} to get $v^1_j, v^2_j \in SH_{\om,m} \cap L^\infty (\Om)$ which are solutions of the equations 
	$$H_{m,\o} (v^1_j)=\varphi^1_j d\mu, \lim\limits_{\Om \ni z\to x} v^1_j (z)=\phi (x)\,\forall x\in\pa \Om$$
and	
	$$H_{m,\o} (v^2_j)=\varphi^2_j d\mu, \lim\limits_{\Om \ni z\to x} v^2_j (z)=\phi (x)\,\forall x\in\pa \Om.$$
According to Theorem \ref{sub}, there exists $\tilde u \in SH_{\om,m}(\Om)\cap L^{\infty}(\Om)$ be such that
$$\lim\limits_{\Om \ni z\to x} \tilde u (z)=\phi (x)\,\forall x\in\pa \Om, \
H_{m,\o}(\tilde u)=F(u,z)d\mu.$$ 	
	Then, using  Theorem \ref{comparison} we see that $v^1_j \downarrow v^1, v^2_j \uparrow v^2$, furthermore, in view of (i) we also have
	\begin{equation} \label{eq2}
		v^1_j \ge u_{j+1} \ge v^2_j.
	\end{equation}
	Next we use (ii) to get 
	$$H_{m,\o} (v^1_j) \to F (u,z) \mu, H_{m,\o} (v^2_j) \to F(u,z)\mu.$$
	So by the monotone convergence theorem (see Theorem 3.4 in \cite{GN18}) we infer 
	$$H_{m} (v^1)=H_m ((v^2)^*)= F(u(z),z)d\mu=H_{m,\o}(\tilde u).$$
	Applying again Theorem \ref{comparison} we obtain $v^1=(v^2)^*=\tilde u$ on $\O.$ \\
Moreover, we have the following easy estimate
\begin{equation} \label{eq1}
|u_{j+1}- \tilde u| \le \max \{v^1_j-\tilde u, \tilde u-v^2_j\}
\end{equation}	
 Fix $\de>0,$ by (\ref{eq1}), for each compact subset $K$ of $\Om$ we have
		$$Cap (K \cap \{|u_{j+1} -\tilde u|>\de\}) \le Cap (K \cap \{v^1_j-\tilde u>\de\})+
	Cap (K \cap \{\tilde u- v^2_j >\de\}).$$
But using Remark \ref{monocap} we see that both monotone sequences $v^1_j$ and $v^2_j$ converge
to $\tilde u$ in $(\o,m)$-capacity. Therefore, we obtain
$$\lim_{j \to \infty} Cap (K \cap \{|u_{j+1} -\tilde u|>\de\})=0.$$
Hence, we have $u_j$ tends to $\tilde u$ in $(\o,m)$-capacity.\\
On the other hand, according to the monotone convergence Theorem, we obtain
$$\lim_{j \to \infty} \int v^1_j d\mu=\lim_{j \to \infty} \int v^2_j d\mu=\int \tilde u d\mu.$$

Notice that, for the second equality, we use Remark \ref{remark} the fact that $\mu$ puts no mass on $m-$polar set.
It follows that
$$\lim_{j \to \infty}\int |u_{j+1}-\tilde u|d\mu=0.$$
Note that, we also have  $u_j \to u$ a.e. $(d\mu).$
Thus, we infer that $\tilde u=u.$ \\
By Proposition \ref{cap} we also get
$$F(u_{j+1},z)d\mu= H_{m,\o} (u_{j+2}) \to H_{m,\o}(\tilde u) \ \text{weakly as}\  j \to \infty.$$
Moreover, according to Lebesgue dominated convergence Theorem, we imply that 
$F(u_j,z) d\mu$ converges weakly to $F(\tilde u,z) d\mu.$
Putting all this together, we imply that $F(\tilde u,z)d\mu=H_m (\tilde u)$.
Therefore, $\tilde u$ is a solution of the problem (*). 
		
\n In the case, $ t \mapsto F(t,z)$ is non-decreasing for every 
$z \in \O \setminus Y$ with $C_m (Y)=0,$ then the uniqueness of $u$ follows directly from Proposition \ref{md4} (with $F=\tilde F, \mu=\nu$).	The proof is complete.
\end{proof}	
	\begin{proof} [{\it The proof of Theorem \ref{th1.2}}]
Let $h\in SH_{\om,m}(\Om)\cap L^{\infty}({\Om})$ be a function that satisfies (\ref{eqh}). Then by inequality \eqref{condi}, for each $j \ge 1$ we have 
$$H_{m,\o} (v+h) \ge H_m (v) \ge Gd\mu \ge F(u_j, z)d\mu=H_{m,\o} (u_j) \ge H_m (h).$$
Thus using Theorem \ref{comparison} we obtain 
$$v+h \le u_j \le h.$$
In particular, the sequence $u_j$ is uniformly bounded in $SH_{\om,m} (\Om) \cap L^\infty (\Om).$
Thus, by the same reasoning as in Theorem \ref{th1.1}, we can suppose, after passing to a subsequence, that
$u_j$ converges pointwsise a.e. $(d\mu)$ to a function $u \in SH_{\om,m} (\Om) \cap L^\infty (\Om).$
Observe that, by the assumption that $F_j (t, \cdot) \to F(t, \cdot)$ locally unniformly on $\R$ when $z \in \Om \setminus Y$ is fixed, we obtain
$$\lim_{j \to \infty} F_j (u_j (z), z)= F(u(z), z), \text{a.e.}\ d \mu.$$
\n
Now define for $z \in \Om,$ the following sequences of non-negative measurable functions
$$\varphi^1_j (z):= \inf_{k \ge j} F_k (u_k (z),z),
\varphi^2_j (z):= \sup_{k \ge j} F_k (u_k (z),z).$$
Then, by inequality \eqref{condi} we obtain:\\
\n 
(i) $0 \le \varphi^1_j (z) \le  F_j (u_j (z),z) \le \varphi^2_j (z) \le G(z)$ for $j \ge 1;$

\n 
(ii) $\lim\limits_{j \to \infty} \varphi^1_j (z)=\lim\limits_{j \to \infty} \varphi^2_j (z)= F(u(z),z)$ a.e. ($d\mu$).

\n 	
Now,  we apply Theorem \ref{sub} to get $v^1_j, v^2_j \in SH_{\om,m} \cap L^\infty (\Om)$ which are solutions of the equations 
$$H_{m,\o} (v^1_j)=\varphi^1_j d\mu, \lim\limits_{\Om \ni z\to x} v^1_j (z)=\phi (x)\,\forall x\in\pa \Om$$
and	
$$H_{m,\o} (v^2_j)=\varphi^2_j d\mu, \lim\limits_{\Om \ni z\to x} v^2_j (z)=\phi (x)\,\forall x\in\pa \Om.$$
Applying Theorem \ref{sub} once again, there exists a function  $\tilde u \in SH_{\om,m}(\Om)\cap L^{\infty}(\Om)$  such that
$$\lim\limits_{\Om \ni z\to x} \tilde u (z)=\phi (x)\,\forall x\in\pa \Om, \
H_{m,\o}(\tilde u)=F(u,z)d\mu.$$ 	
Then, repeating the argument as in the proof of Theorem \ref{th1.1} we infer that 
$v^1_j \downarrow v^1, v^2_j \uparrow v^2$, and
$v^1=(v^2)^*=\tilde u$ on $\O.$  Moreover, 
$u_j \to \tilde u$ in $(\o,m)$-capacity on $\Om,$ and $u=\tilde u$ a.e. $(d\mu).$
Finally, we apply Proposition \ref{cap} and the Lebesgue dominated convergence theorem to see that 
$F(\tilde u,z)d\mu=H_{m,\o} (\tilde u)$.
Thus, $\tilde u$ is indeed, a solution of the problem (*). The proof is complete.
\end{proof}		

	\section*{Declarations}
	\subsection*{Ethical Approval}
	This declaration is not applicable.
	\subsection*{Competing interests}
	The authors have no conflicts of interest to declare that are relevant to the content of this article.
	\subsection*{Authors' contributions }
Hoang Thieu Anh,	Le Mau Hai, Nguyen Quang Dieu and Nguyen Van Phu  together studied  the manuscript.
	\subsection*{Availability of data and materials}
	This declaration is not applicable.

\end{document}